\documentclass[12pt]{article}
\usepackage{amsthm, amssymb, amsmath, latexsym}
\usepackage[utf8]{inputenc}
\usepackage[english]{babel}
\usepackage{longtable}

\begin{document}
\pagestyle{myheadings}
\thispagestyle{empty}
\setcounter{page}{1}

\newtheorem{definition}{Definition}
\newtheorem{proposition}{Proposition}
\newtheorem{theorem}{Theorem}
\newtheorem{lemma}{Lemma}
\newtheorem{corollary}{Corollary}
\newtheorem{remark}{Remark}
\theoremstyle{plain}
\mathsurround 2pt

\gdef\Aut{\mathop{\rm Aut}\nolimits}
\gdef\Ker{\mathop{\rm Ker}\nolimits}
\gdef\Im{\mathop{\rm Im}\nolimits}
\gdef\End{\mathop{\rm End}\nolimits}
\gdef\Inn{\mathop{\rm Inn}\nolimits}
\gdef\exp{\mathop{\rm exp}\nolimits}

\title{Lower bounds for the number of local nearrings on groups of order $p^3$}

\author{Iryna Raievska, Maryna Raievska\\
University of Warsaw, Warsaw, Poland;\\
Institute of Mathematics of NAS of Ukraine, Ukraine\\
raeirina@imath.kiev.ua, raemarina@imath.kiev.ua}

\date{}

\maketitle

\begin{abstract}
Lower bounds for the number of local nearrings on groups of order $p^3$ are obtained. On each non-metacyclic non-abelian or metacyclic abelian groups of order $p^3$ there exist at least $p+1$ non-isomorphic local nearrings.
\end{abstract}

\section{Introduction}

\

A study of local nearrings was first initiated in~\cite{MCJ_68} and it was found that the additive group of a finite zero-symmetric local nearring is a $p$-group. In~\cite{CM_68} it is shown that, up to an isomorphism, there exist $p-1$ local zero-symmetric nearrings with elementary abelian additive groups of order $p^2$, in which the subgroups of non-invertible elements have order $p$, that is, those nearrings which are not nearfields. Together with the fundamental paper~\cite{Z_36} and \cite{ClMal_66}, it is obtained a complete description of all zero-symmetric local nearrings of order $p^2$. For instance, every nearring with identity on a cyclic group is a commutative ring.

Note that there is no nearring with identity whose additive group is isomorphic to the quaternion group $Q_8$~\cite{Clay_70}. The dihedral group $D_4$ of order $8$ cannot be the additive group of local nearrings~\cite{CM_71}. The existence of local nearrings on finite abelian $p$-groups is proved in~\cite{CM_70}, i.e. every non-cyclic abelian $p$-group of order $p^n>4$ is the additive group of a zero-symmetric local nearring which is not a ring. Also, it is established in \cite{IMYa_12} that an arbitrary non-metacyclic Miller--Moreno $p$-group of order $p^n>8$ is the additive group of some local nearring, and the multiplicative group of such nearring has order $p^{n-1}(p-1)$. All nearrings with identity up to the order of $31$ are contained in the package SONATA~\cite{SONATA} of the computer algebra system~GAP~\cite{GAP4}.

In \cite{IM_2021} it is proved that, up to an isomorphism, there exist at least $p$ local nearrings on elementary abelian additive groups of order $p^3$, which are not nearfields. Lower bounds for the number of local nearrings on groups of order $p^3$ are obtained. It is established that on each non-metacyclic non-abelian or metacyclic abelian groups of order $p^3$ there exist at least $p+1$ non-isomorphic local nearrings.

\section{Preliminaries}

We will give the basic definitions.

\begin{definition}
A non-empty set $R$ with two binary operations $``+"$ and $``\cdot"$ is a \textbf{nearring} if:
\begin{description}
  \item[1)] $(R,+)$ is a group with neutral element $0${\rm ;}
  \item[2)] $(R,\cdot)$ is a semigroup{\rm ;}
  \item[3)] $x\cdot (y+z)=x\cdot y+x\cdot z$ for all $x$, $y$, $z\in R$.
\end{description}
Such a nearring is called a left nearring. If axiom 3) is replaced by an axiom $(x+y)\cdot z = x\cdot z + y\cdot z$ for all $x$, $y$, $z\in R$, then we get a right
nearring.
\label{nr}
\end{definition}

The group $(R,+)$ of a nearring $R$ is denoted by $R^+$ and called the {\em additive group} of $R$. It is easy to see that for each subgroup $M$ of $R^+$ and for each element $x\in R$ the set $xM=\{x\cdot y|y\in M\}$ is a subgroup of $R^+$ and in particular $x\cdot 0=0$. If in addition $0\cdot x=0$ for all $x\in R$, then the nearring $R$ is called {\em zero-symmetric}. Furthermore, $R$ is a {\em nearring with an identity} $i$  if the semigroup $(R,\cdot)$ is a monoid with identity element $i$. In the latter case the group of all invertible elements of the monoid $(R,\cdot)$ is denoted by $R^*$ and called the {\em multiplicative group} of $R$.  A subgroup $M$ of $R^+$ is called $R^*$-{\em invariant}, if $rM\leq M$ for each $r\in R^*$, and $(R,R)$-{\em subgroup}, if $xMy\subseteq M$ for arbitrary $x$, $y\in R$.\medskip

The following assertion is well-known (see, for instance, \cite{ClMal_66}, Theorem 3).

\begin{lemma}
The exponent of the additive group of a finite nearring $R$ with identity $i$ is equal to the additive order of $i$ which coincides with the additive order of every invertible element of $R$.
\label{exponent}
\end{lemma}

\begin{definition}
A nearring $R$ with identity is called \textbf{local} if the set $L$ of all non-invertible elements of $R$ forms a subgroup of the additive group $R^{+}$.
\end{definition}

Through this paper $L$ will denote the subgroup of non-invertible elements of $R$.

The following lemma characterizes the main properties of finite local nearrings (see~\cite{AHS_2004}, Lemma~3.2).

\begin{lemma}
Let $R$ be a local nearring with identity $i$. Then the following statements hold{\rm :}
\begin{description}
\item[1)] $L$ is an $(R,R)$-subgroup of $R^{+}${\rm ;}
\item[2)] each proper $R^*$-invariant subgroup of $R^+$ is contained in $L${\rm ;}
\item[3)] the set $i+L$ forms a subgroup of the multiplicative group $R^*$.
\end{description}
\label{prop}
\end{lemma}

Finite local nearrings with a cyclic subgroup of non-invertible elements are described in~\cite[Theorem~1]{RIM_2}.

\begin{theorem}
Let $R$ be a local nearring of order $p^n$ with {$n\!> 1$} whose subgroup $L$ is cyclic and non-trivial. Then the additive group $R^+$ is either cyclic or is an elementary abelian group of order $p^2$. In the first case, $R$ is a commutative local ring, which is isomorphic to residual ring $\mathbb Z/p^n\mathbb Z$ with $n\ge 2$, in the other case there exist $p$ non-isomorphic such nearrings $R$ with $|L|=p$, from which $p-1$ are zero-symmetric nearrings and their multiplicative groups $R^{*}$ are isomorphic to a semidirect product of two cyclic subgroups of orders $p$ and $p-1$.
\label{theorem_2}
\end{theorem}

As a direct consequence of Theorem~\ref{theorem_2} we have the following result.

\begin{corollary}
Let $R$ be a local nearring of order $p^3$ which is not isomorphic to $\mathbb Z/p^3 \mathbb Z$ or is not a nearfield. Then the subgroup of non-invertible elements $L$ is an elementary abelian group of order $p^2$.
\label{cor_1}
\end{corollary}

The following statement contains a classification of groups of order $p^3$ (see~\cite{Hall_1962}).

\begin{proposition}
Let $G$ be a group of order $p^3$. The defining relations of such non-isomorphic groups are given:

Abelian groups{\rm :}

\begin{description}
  \item[1)] ${a^p}^3=1$.
  \item[2)] ${a^p}^2=1$, $b^p=1$, $ab=ba$.
  \item[3)] $a^p=b^p=c^p=1$, $ab=ba$, $ac=ca$, $cb=bc$.
\end{description}

Non-abelian groups of order $2^3=8${\rm :}

\begin{description}
  \item[4)] a dihedral group, $a^4=1$, $b^2=1$, $a^{-1}b=ba$.
  \item[5)] a quaternion group, $a^4=1$, $b^2=a^2$, $a^{-1}b=ba$.
\end{description}

Non-abelian groups of order $p^3$, $p$ is odd{\rm :}

\begin{description}
  \item[6)] ${a^p}^2=1$, $b^p=1$, $b^{-1}ab=a^{1+p}$.
  \item[7)] $a^p=1$, $b^p=1$, $c^p=1$, $ab=bac$, $ac=ca$, $bc=cb$.
\end{description}
\label{rem1}
\end{proposition}

Next, we denote by $G_1$ a group with relations~7), $G_2$ a group with relations~6), $G_3$ a group with relations~2) of Proposition~\ref{rem1}.

The following two lemmas are appeared in \cite{IMYa_12}.

\begin{lemma}
Let $G_1$ be an additively written group. Then for any natural numbers $k$ and $l$ the equalities $-ak-bl+ak+bl=c(kl)$ and $bl+ak=-c(kl)+ak+bl$ hold.
\label{4.1c}
\end{lemma}

\begin{proof}
Since $-b+a+b=a+c$, we have $-bl+a+bl=a+cl$. Then
$$-bl+ak+bl=(a+cl)k=ak+ckl.$$
Therefore, $-ak-bl+ak+bl=ckl$.
\end{proof}

\begin{lemma}
Let $G_1$ be an additively written group. Then for any natural numbers $k$, $l$ and $r$ the equality $(ak+bl)r=akr+blr-ckl\binom{r}{2}$ holds.
\label{4.3c}
\end{lemma}

\begin{proof}
The proof will be carried out by induction on $r$. For $r=1$ the equality is valid.
Let for $r$ the equality hold, i.e.
$$(ak+bl)r=akr+blr-ckl\binom{r}{2}.$$
Let us prove the equality for $r+1$:
$$\begin{array}{c}
(ak+bl)(r+1)=akr+blr+ak+bl-ckl\binom{r}{2}=\\
=ak(r+1)+bl(r+1)-cklr-ckl\binom{r}{2}=\\
=ak(r+1)+bl(r+1)-ckl(r+\binom{r}{2})=\\
=ak(r+1)+bl(r+1)-ckl\binom{r+1}{2}.
\end{array}$$
Therefore, the equality is valid for any $r$.
\end{proof}

Let additively written groups of type $H$ have a finite representation in the form
$$\langle a, b| ap^k, bp, -b+a+b-a(1+p^{k-1})\rangle,$$
where $k\geq 2$ and $p$ is prime (see~\cite{LL_80}).

The number of non-isomorphic nearrings with identity on groups of type $H$ is given in~\cite{LL_80}. It is obvious that for $k=2$ and $p>2$ a group of type $H$ will be isomorphic to the group $G_2$.

As noted above, there exist local nearrings on all abelian groups. Also, according to~\cite{IMYa_12} and \cite{LL_80}, there exist local nearrings on $G_1$ and $G_2$, respectively. So, we have the following result.

\begin{proposition}
On each group of order $p^3$ with $p>2$ there exists a local nearring.
\label{prop2}
\end{proposition}

Denote by $n(G)$ the number of all non-isomorphic local nearrings on the group $G$.

\section{Nearrings with identity whose additive groups are isomorphic to $G_1$}

Let $R$ be a nearring with identity whose additive group of $R^+$ is isomorphic to $G_1$. Then $R^{+}=\langle a\rangle +\langle b \rangle +\langle c \rangle$ for some elements $a$, $b$ and $c$ of $R$ satisfying the relations $ap=0$, $bp=0$, $cp=0$, $a+b=b+a+c$, $a+c=c+a$ and $b+c=c+b$. In particular, each element $x\in R$ is uniquely written in the form $x=ax_1+bx_2+cx_3$ with coefficients $0\le x_1<p$, $0\le x_2<p$ and $0\le x_3<p$.

Since the order of the element $a$ is equal to the exponent of group $G$, then by Lemma~\ref{exponent} we can assume that $a$ is an identity of $R$, i.e. $ax=xa=x$ for each $x\in R$. Furthermore, for each $x\in R$ there exist coefficients $\alpha(x)$, $\beta(x)$ and $\gamma(x)$ such that $xb=a\alpha(x)+b\beta(x)+c\gamma(x)$. It is clear that they are uniquely defined modulo $p$, so that some mappings $\alpha\colon R\to \mathbb Z_{p}$, ${\beta: R\rightarrow Z_{p}}$ and ${\gamma: R\rightarrow Z_p}$ are determined.

Nearrings with identity and local nearrings on non-metacyclic Miller-Moreno groups were studied in \cite{IMYa_12} and \cite{IMYa_11}.  Lemmas \ref{idN}, \ref{G_1_id}, \ref{G_3_id} are based on the results of these papers.

\begin{lemma}
Let $R$ be a nearring with identity whose additive group of $R^+$ is isomorphic to $G_1$.  If $a$ coincides with identity element of $R$, $x=ax_1+bx_2+cx_3$, $y=ay_1+by_2+cy_3\in R$, $xb=a\alpha(x)+b\beta(x)+c\gamma(x)$ then
$$xy=a(x_1y_1+\alpha(x)y_2)+b(x_2y_1+\beta(x)y_2)+c(-x_1x_2\binom{y_1}{2}-$$
$$-\alpha(x)\beta(x)\binom{y_2}{2}-x_2\alpha(x)y_1y_2+x_3y_1+\gamma(x)y_2+x_1\beta(x)y_3-x_2\alpha(x)y_3).$$
Moreover, for the mappings $\alpha\colon R\to \mathbb Z_{p}$, $\beta\colon R\to \mathbb Z_{p}$ and $\gamma\colon R\to \mathbb Z_{p}$ the following statements hold{\rm :}
\begin{description}
  \item[(0)] $\alpha(0)\equiv 0\; (\!\!\mod p)$, $\beta(0)\equiv 0\; (\!\!\mod p)$ and ${\gamma(0)\equiv 0\; (\!\!\mod p)}$ if and only if the nearring $R$ is zero-symmetric;
  \item[(1)] $\alpha(xy)\equiv x_1\alpha(y)+\alpha(x)\beta(y)\; (\!\!\mod p\;),$
  \item[(2)] $\beta(xy)\equiv x_2\alpha(y)+\beta(x)\beta(y)\; (\!\!\mod p\;),$
  \item[(3)] $\gamma(xy)\equiv -x_1x_2\binom{\alpha(y)}{2}-\alpha(x)\beta(x)\binom{\beta(y)}{2}-x_2\alpha(x)\alpha(y)\beta(y)+$
$$+x_3\alpha(y)+\gamma(x)\beta(y)+x_1\beta(x)\gamma(y)-x_2\alpha(x)\gamma(y)\; (\!\!\mod p\;).$$
\end{description}
\label{idN}
\end{lemma}

\begin{proof}
Since $0\cdot a=a\cdot 0=0$, it follows that $R$ is a zero-symmetric nearring if and only if $$0=0\cdot b=a\alpha(0)+b\beta(0)+c\gamma(0)$$ or equivalently $\alpha(0)\equiv 0\; (\!\!\mod p)$, ${\beta(0)\equiv 0\; (\!\!\mod p)}$ and ${\gamma(0)\equiv 0\; (\!\!\mod p)}$. Moreover, since $c=-a-b+a+b$ and the left distributive law we have $0\cdot c=-0\cdot a-0\cdot b+0\cdot a+0\cdot b=0$, whence $$0\cdot x=0\cdot (ax_1+bx_2+cx_3)=(0\cdot a)x_1+(0\cdot b)x_2+(0\cdot c)x_3=0.$$ So that statement (0) holds.

Further, using Lemma~\ref{4.1c}, we derive

$$\begin{array}{c}
xc=-xa-xb+xa+xb=-cx_3-bx_2-ax_1-c\gamma(x)-b\beta(x)-a\alpha(x)+\\
+ax_1+bx_2+cx_3+a\alpha(x)+b\beta(x)+c\gamma(x)=\\
=-bx_2-ax_1-b\beta(x)-a\alpha(x)+ax_1+bx_2+a\alpha(x)+b\beta(x)=\\
=-bx_2+cx_1\beta(x)-b\beta(x)-ax_1-a(\alpha(x)-x_1)+bx_2+a\alpha(x)+b\beta(x)=\\
=cx_1\beta(x)-b(x_2+\beta(x))-a\alpha(x)+bx_2+a\alpha(x)+b\beta(x)=\\
=cx_1\beta(x)-b(x_2+\beta(x))-a\alpha(x)-cx_2\alpha(x)+a\alpha(x)+bx_2+b\beta(x)=\\
=c(x_1\beta(x)-x_2\alpha(x))-b(x_2+\beta(x))+bx_2+b\beta(x)=c(x_1\beta(x)-x_2\alpha(x)).
\end{array}$$

Further, using the left distributive law, we obtain
$$xy=(ax_1+bx_2+cx_3)y_1+(a\alpha(x)+b\beta(x)+c\gamma(x))y_2+$$
$$+(cx_1\beta(x)-x_2\alpha(x))y_3.$$
By Lemma~\ref{4.3c}, we get
$$(ax_1+bx_2)y_1=ax_1y_1+bx_2y_1-cx_1x_2\binom{y_1}{2},$$
$$(a\alpha(x)+b\beta(x))y_2=a\alpha(x)y_2+b\beta(x)y_2-c\alpha(x)\beta(x)\binom{y_2}{2}$$
and
$$bx_2y_1+a\alpha(x)y_2=a\alpha(x)y_2+bx_2y_1-cx_2\alpha(x)y_1y_2.$$
Hence and using the left distributive law, we have
$$xy=a(x_1y_1+\alpha(x)y_2)+b(x_2y_1+\beta(x)y_2)+c(-x_1x_2\binom{y_1}{2}-$$
$$-\alpha(x)\beta(x)\binom{y_2}{2}-x_2\alpha(x)y_1y_2+x_3y_1+\gamma(x)y_2+x_1\beta(x)y_3-x_2\alpha(x)y_3).$$
The associativity of multiplication in $R$ implies that for all $x$, $y\in R$ $$(xy)b=x(yb).\leqno 1)$$

According to $xb=a\alpha(x)+b\beta(x)+c\gamma(x)$, we obtain $$(xy)b=a\alpha(xy)+b\beta(xy)+c\gamma(xy)\leqno 2)$$ and
$yb=a\alpha(y)+b\beta(y)+c\gamma(y)$. Substituting the last equation to the right part of equality 1), we also have
$$x(yb)=a(x_1\alpha(y)+\alpha(x)\beta(y))+b(x_2\alpha(y)+\beta(x)\beta(y))+\leqno 3)$$
$$+c(-x_1x_2\binom{\alpha(y)}{2}-\alpha(x)\beta(x)\binom{\beta(y)}{2}-x_2\alpha(x)\alpha(y)\beta(y)+$$
$$+x_3\alpha(y)+\gamma(x)\beta(y)+x_1\beta(x)\gamma(y)-x_2\alpha(x)\gamma(y)).$$
Since equality 1) implies the congruence of the corresponding coefficients in formulas 2) and 3), we obtain statements (1)--(3).
\end{proof}

\section{Local nearrings whose additive groups are isomorphic to $G_1$}

Let $R$ be a local nearring whose additive group of $R^+$ is isomorphic to $G_1$. Then $R^{+}=\langle a\rangle +\langle b \rangle +\langle c \rangle$ for some elements $a$, $b$ and $c$ of $R$ satisfying the relations $ap=0$, $bp=0$, $cp=0$, $a+b=b+a+c$, $a+c=c+a$ and $b+c=c+b$. In particular, each element $x\in R$ is uniquely written in the form $x=ax_1+bx_2+cx_3$ with coefficients $0\le x_1<p$, $0\le x_2<p$ and $0\le x_3<p$.

Since order of the element $a$ is equal to the exponent of group $G$, then be Lemma~\ref{exponent} we can assume that $a$ is an identity of $R$, i.~e. $ax=xa=x$ for each $x\in R$. Furthermore, for each $x\in R$ there exist coefficients $\alpha(x)$, $\beta(x)$ and $\gamma(x)$ such that $xb=a\alpha(x)+b\beta(x)+c\gamma(x)$. It is clear that they are uniquely defined modulo $p$, so that some mappings ${\alpha: R\rightarrow Z_{p}}$, ${\beta: R\rightarrow Z_{p}}$ and ${\gamma: R\rightarrow Z_p}$ are determined.

By Corollary~\ref{cor_1}, $L$ is the normal subgroup of order $p^2$ in $R$. Since $L$ consists the derived subgroup of $R^+$ it follows that the generators $b$ and $c$ we can choose such that $c=-a-b+a+b$. Then $L=\langle b\rangle + \langle c \rangle$ and the subgroup $\langle c \rangle$ is the center of $R^+$. Since $R^*=R\setminus L$  it follows ${R^*=\{ax_1+bx_2+cx_3\mid x_1\not\equiv 0 \; (\!\!\mod p\;)\}}$ and $x=ax_1+bx_2+cx_3$ is invertible if and only if $x_1\not\equiv 0 \; (\!\!\mod p\;)$.

\begin{lemma}
Let $R$ be a local nearring whose additive group of $R^+$ is isomorphic to $G_1$. If $a$ coincides with identity element of $R$, $x=ax_1+bx_2+cx_3$, $y=ay_1+by_2+cy_3\in R$, $xb=a\alpha(x)+b\beta(x)+c\gamma(x)$, then
$$x\cdot y=ax_1y_1+b(x_2y_1+\beta(x)y_2)+$$
$$+c(-x_1x_2\binom{y_1}{2}+x_3y_1+\gamma(x)y_2+x_1\beta(x)y_3).~\rm{(*)}$$
Moreover, for the mappings $\beta\colon R\to \mathbb Z_{p}$ and $\gamma\colon R\to \mathbb Z_{p}$ the following statements hold{\rm :}
\begin{description}
 \item[\rm{(0)}] $\alpha(0)\equiv 0\; (\!\!\mod p)$, $\beta(0)\equiv 0\; (\!\!\mod p)$ and ${\gamma(0)\equiv 0\; (\!\!\mod p)}$ if and only if the nearring $R$ is zero-symmetric{\rm ;}
 \item[\rm{(1)}] $\alpha(x)\equiv 0\; (\!\!\mod p)${\rm ;}
 \item[\rm{(2)}] if $\beta(x)\equiv 0 \; (\!\!\mod p)$, then $x_1\equiv 0 \; (\!\!\mod p)${\rm ;}
 \item[\rm{(3)}] $\beta(xy)\equiv \beta(x)\beta(y)\; (\!\!\mod p\;)${\rm ;}
 \item[\rm{(4)}] $\gamma(xy)\equiv \gamma(x)\beta(y)+x_1\beta(x)\gamma(y)\; (\!\!\mod p)$.
\end{description}
\label{G_1_id}
\end{lemma}

\begin{proof}
Since $L=\langle b\rangle + \langle c \rangle$ and $L$ is the $(R,R)$-subgroup in $R^+$ by statement~1) of Lemma~\ref{prop} it follows that $xb\in L$, hence $a\alpha(x)\in L$ for each $x\in R$. Thus $\alpha(x)\equiv 0 \; (\!\!\mod p)$ and we get statement~(1). Substituting the obtained value of $\alpha(x)\equiv 0 \; (\!\!\mod p)$ in the formulas from Lemma~\ref{idN}, we obtain statement~(3) and (4) of the lemma and the formula for the product $xy$. Putting $y=c$, we get $xc=c(x_1\beta(x))$. Hence, if $\beta(x)\equiv 0 \; (\!\!\mod p)$, then $xc=0$, and so $x\in L$.  Therefore, $x_1\equiv 0 \; (\!\!\mod p)$, as claimed in statement (2). Indeed, statement (0) repeats the statement (0) of Lemma~\ref{idN}.
\end{proof}

It is known that for such groups the commutator $D(R^+)$ coincides with the center $Z(R^+)$ and has the order $p$.

\begin{lemma}
The commutator $D(R^+)$ is an ideal in the local nearring $R$.
\label{prop3}
\end{lemma}

\begin{proof}
Let $x=ax_1+bx_2+cx_3$, $y=ay_1+by_2+cy_3\in R$, $z=cz_3\in D(R^+)$. Let us check whether $D(R^+)$ is an ideal in $R$, i.~e. $(z+x)y-xy\in D(R^+)$. To do this, we use formula $\rm{(*)}$ for multiplying elements in $R$.
We obtain
$$\begin{array}{c}
   (z+x)y-xy=(cz_3+ax_1+bx_2+cx_3)(ay_1+by_2+cy_3)-\\
   (ax_1+bx_2+cx_3)(ay_1+by_2+cy_3)=(ax_1+bx_2+c(x_3+z_3))(ay_1+by_2+cy_3)-\\
   (ax_1+bx_2+cx_3)(ay_1+by_2+cy_3)=ax_1y_1+\\
   b(x_2y_1+\beta(x)y_2)+c(-x_1x_2\binom{y_1}{2}+x_3y_1+z_3y_1+\gamma(z+x)y_2+x_1\beta(x)y_3)-\\
   ax_1y_1-b(x_2y_1+\beta(x)y_2)-c(-x_1x_2\binom{y_1}{2}+x_3y_1+\gamma(x)y_2+x_1\beta(x)y_3)=\\
   ax_1y_1+b(x_2y_1+\beta(x)y_2)-ax_1y_1-b(x_2y_1+\beta(x)y_2)+c(z_3y_1+(\gamma(z+x)-\gamma(x))y_2)=\\
   ax_1y_1-ax_1y_1+b(x_2y_1+\beta(x)y_2)-c(x_1y_1(x_2y_1+\beta(x)y_2))-\\
   b(x_2y_1+\beta(x)y_2)+c(z_3y_1+(\gamma(z+x)-\gamma(x))y_2)=\\
   c(z_3y_1+(\gamma(z+x)-\gamma(x))y_2-x_1x_2y_1^2-x_1y_1y_2\beta(x))\in D(R^+).
\end{array}$$
Therefore, $D(R^+)$ is an ideal of $R$.
\end{proof}

\begin{lemma}
Let $R$ be a local nearring whose additive group of $R^+$ is isomorphic to $G_1$. If $x=ax_1+bx_2+cx_3$, $y=ay_1+by_2+cy_3\in R$, then the mappings $\beta\colon R\to \mathbb Z_{p^2}$ and $\gamma\colon R\to \mathbb Z_{p}$ from $\rm{(*)}$ can be one of the following{\rm :}
\begin{description}
  \item[1)] $\beta(x)={x_1}^i$ and $\gamma(x)=0$ ($0<i<p$){\rm ;}
  \item[2)] $\beta(x)=1$ and $\gamma(x)=0${\rm ;}
  \item[3)] $\beta(x)=x_1^2$ and $\gamma(x)=x_1x_2$.
\end{description}
\label{lemma3}
\end{lemma}

\begin{proof}
1) Since zero-symmetric local nearrings of order $p^2$ are classified in~\cite{CM_68} it follows that all non-isomorphic factor-nearrings with derived subgroup $N=R/D(R^+)$ are described. That is, you can apply the multiplication formula from the specified work, pre-adapting it for the left local nearrings. It is clear that $N^+=\langle a \rangle + \langle b \rangle$. Namely, let $x=ax_1+bx_2$ and $y=ay_1+ by_2$ be elements of $N$, then

$$xy=ax_1y_1+b(x_2y_1+\rho(x_1)y_2).$$

Moreover, by \cite[Theorem~1.6]{CM_68} $\rho$ takes one of $p-1$ values for zero-symmetric nearrings.

On the other hand, by the formula~$\rm{(*)}$ we have:

$$xy=ax_1y_1+b(x_2y_1+\beta(x)y_2).$$

Equating the coefficients for generators, we obtain: $\beta(x)y_2=\rho(x_1)y_2$. Hence $\beta(x)=\rho(x_1)$. So there exist $p-1$ different zero-symmetric nearrings.

2) It is obvious that the multiplication~$\rm{(*)}$ with the functions $\beta(x)=1$ and $\gamma(x)=0$ is the constant nearring multiplication.

3) We show further that for $\beta(x) = x_1^2$ and $\gamma(x) = x_1x_2$ the multiplication~$\rm{(*)}$ is a nearring multiplication. By Lemma~\ref{G_1_id} (3) and (4), the functions $\beta(x)$ and $\gamma(x)$ satisfy the following conditions: $\beta(xy)=\beta(x)\beta(y)$ and $\gamma(xy)=\beta(y)\gamma(x)+x_1\beta(x)\gamma(y)$. Hence $\beta(xy)=x_1^2y_1^2=\beta(x)\beta(y)$ and $\gamma(xy)=x_1y_1(x_2y_1+x_1y_2)=x_1x_2y_1^2+x_1^2y_1y_2=\beta(y)\gamma(x)+x_1\beta(x)\gamma(y)$.
\end{proof}

So, we have examples of $p+1$ nearring multiplications and we formulate the following result.

\begin{theorem}
There exist at least $p+1$ non-isomorphic local nearrings on each non-metacyclic non-abelian groups of order $p^3$.
\label{theorem_3}
\end{theorem}

\begin{proof}
The nearrings $R=(R,+,\cdot)$ with functions $\beta(x)={x_1}^i$ and $\gamma(x)=0$ ($0<i<p$) are non-isomorphic zero-symmetric local nearrings according to the above and the paper~\cite{CM_68}. It is easy to check that the local nearrings $R=(R,+,\cdot)$ with functions $\beta(x)=1$ and $\gamma(x)=0$ are non-zero-symmetric. It is obvious that $R=(R,+,\cdot)$ with functions $\beta(x)=x_1^2$ and $\gamma(x)=x_1x_2$ are zero-symmetric and non-isomorphic to the nearrings considered above. Therefore, there exist at least $p+1$  non-isomorphic local nearrings on each non-metacyclic non-abelian groups of order $p^3$.
\end{proof}

\textbf{Example 1.} Let $G\cong (C_5\times C_5)\rtimes C_5$. If $x=ax_1+bx_2+cx_3$ and $y=ay_1+by_2+cy_3\in G$ and $(G,+, \cdot)$ is a local nearring, then by Lemma~\ref{lemma3} $``\cdot"$ is one of the following multiplications.
\begin{description}
  \item[(1)] $x\cdot y=ax_1y_1+b(x_2y_1+y_2)+c(-x_1x_2\binom{y_1}{2}+x_3y_1+x_1y_3);$
  \item[(2)] $x\cdot y=ax_1y_1+b(x_2y_1+x_1^3(x)y_2)+c(-x_1x_2\binom{y_1}{2}+x_3y_1+x_1x_2y_2+x_1^4y_3);$
  \item[(3)] $x\cdot y=ax_1y_1+b(x_2y_1+x_1y_2)+c(-x_1x_2\binom{y_1}{2}+x_3y_1+x_1^2y_3);$
  \item[(4)] $x\cdot y=ax_1y_1+b(x_2y_1+x_1^2y_2)+c(-x_1x_2\binom{y_1}{2}+x_3y_1+x_1^3y_3);$
  \item[(5)] $x\cdot y=ax_1y_1+b(x_2y_1+x_1^3(x)y_2)+c(-x_1x_2\binom{y_1}{2}+x_3y_1+x_1^4y_3);$
  \item[(6)] $x\cdot y=ax_1y_1+b(x_2y_1+x_1^4y_2)+c(-x_1x_2\binom{y_1}{2}+x_3y_1+x_1^5y_3),$
\end{description}

A computer program verified that for $p=5$, the nearring obtained in Lemma~\ref{lemma3}~(3) is indeed a local nearring (see Example~1~(2)), is deposited on GitHub:

\verb+https://github.com/raemarina/Examples/blob/main/LNR_125-3.txt+

From the packages SONATA and LocalNR~\cite{RRS_21} we have the following number of non-isomorphic local nearrings.

\begin{table}[h]
\centering
\begin{tabular}{|c|c|}
\hline
$StructureDescription(R^+)$ & $n(R^+)$\\
 \hline
$(C_3\times C_3)\rtimes C_3$                     &  4\\
 \hline
$(C_5\times C_5)\rtimes C_5$                     &  6\\
 \hline
$(C_7\times C_7)\rtimes C_7$                     &  8\\
 \hline
$(C_{11}\times C_{11})\rtimes C_{11}$            &  12\\
   \hline
\end{tabular}
\end{table}

\section{Local nearrings whose additive groups are isomorphic to $G_2$}

\

Let $R$ be a local nearring whose additive group of $R^+$ is isomorphic to $G_2$. Then $R^{+}=\langle a\rangle +\langle b \rangle$ for some elements $a$ and $b$ of $R$ satisfying the relations ${a^p}^2=1$, $bp=0$ and $b^{-1}ab=a^{1+p}$. In particular, each element $x\in R$ is uniquely written in the form $x=ax_1+bx_2$ with coefficients $0\le x_1<p^2$ and $0\le x_2<p$.

By~\cite[Theorem~7.1]{LL_80} for $p = 3$ there exist three zero-symmetric nearrings with identity on $G_2$, and for $p> 3$ one. At the same time, by~\cite[Theorem~4.2]{LL_80} there exists one non-zero-symmetric nearring with identity on $G_2$. On the other hand, in~\cite{CM_71} it was shown that for each group $G_2$ there exists a zero-symmetric local nearring. With the SONATA nearring library, it is easy to make sure that all nearrings with identity on $G_2$ of order $27$ are local. The formula for multiplying elements of local nearrings on Miller--Moreno metacyclic groups is defined in~\cite{IYa_12}. Since $G_2$ is a Miller--Moreno metacyclic group, using~\cite[Corollary~2]{IYa_12}, for arbitrary elements $x=ax_1+bx_2$ and $y=ay_1+by_2$ of $G_2$ and putting $\alpha(x)=0$ and $\beta(x)=1$, we obtain the following multiplication formula:
$$x\cdot y=a(x_1y_1-x_1x_2\binom{y_1}{2}p)+b(x_2y_1+y_2).~\rm{(*\;*)}$$
It is easy to see that $R = (G_2,+,\cdot)$ is a non-zero-symmetric local nearring.

The multiplication formula for arbitrary elements of a zero-symmetric local nearring  on $G_2$ is given in the proving of~\cite[Theorem~2]{IYa_12}, namely:
$$x\cdot y=a(x_1y_1-x_1x_2\binom{y_1}{2}p)+b(x_2y_1+\beta(x)y_2),~\rm{(***)}$$
where $\beta(x)=\left\{
               \begin{array}{ll}
                 1, & if~\hbox{$x_1\not\equiv 0\; (\!\!\mod p\;)$;}\\
                 0, & if~\hbox{$x_1\equiv 0\; (\!\!\mod p\;)$.}
               \end{array}
             \right.$

So, from the results presented in this section, we have:
\begin{description}
  \item[1)] $n(G_2)=4$ for $p=3${\rm ;}
  \item[2)] $n(G_2)=2$ for $p>3$. Moreover, the multiplication formulas in such nearrings are determined by formulas~$\rm{(*\;*)}$ or $\rm{(***)}$.
\end{description}

\section{Local nearrings whose additive groups are isomorphic to $G_3$}

Let $R$ be a local nearring, which is not nearfield, whose additive group of $R^+$ is isomorphic to $G_3$. Then $R^{+}=\langle a\rangle +\langle b \rangle $ for some elements $a$ and $b$ of $R$ satisfying the relations $ap^2=0$, $bp=0$ and $b+a=a+b$. In particular, each element $x\in R$ is uniquely written in the form $x=ax_1+bx_2$ with coefficients $0\le x_1<p^2$ and $0\le x_2<p$.

We can without loss of generality assume that $a$ is an identity of $R$, i.~e. $ax=xa=x$ for each $x\in R$. Furthermore, for each $x\in R$ there exist coefficients $\alpha(x)$ and $\beta(x)$ such that $xb=a\alpha(x)+b\beta(x)$. It is clear that they are uniquely defined modulo $p^2$ and $p$, respectively, so that some mappings ${\alpha: R\rightarrow Z_{p^2}}$ and ${\beta: R\rightarrow Z_p}$ are determined.

\begin{lemma}
Let $R$ be a local nearring, which is not nearfield, whose additive group of $R^+$ is isomorphic to $G_3$. If $a$ coincides with identity element of $R$, $x=ax_1+bx_2$, $y=ay_1+by_2\in R$, $xb=a\alpha(x)+b\beta(x)$, then
$$xy=a(x_1y_1+\alpha(x)y_2)+b(x_2y_1+\beta(x)y_2).~\rm{(***\;*)}$$
Moreover, for the mappings $\alpha\colon R\to \mathbb Z_{p^2}$ and $\beta\colon R\to \mathbb Z_{p}$ the following statements hold{\rm :}
\begin{itemize}
 \item[\rm{(0)}] $\alpha(0)=\beta(0)= 0$ if and only if $R$ is zero-symmetric{\rm ;}
 \item[\rm{(1)}] $\alpha(a)= 0$ and $\beta(a)= 1${\rm ;}
 \item[\rm{(2)}] $\alpha(x)\equiv 0\; (\!\!\mod p)${\rm ;}
 \item[\rm{(3)}] $\alpha(xy)\equiv x_1\alpha(y)+\alpha(x)\beta(y)\; (\!\!\mod p\;)${\rm ;}
 \item[\rm{(4)}] $\beta(xy)\equiv x_2\alpha(y)+\beta(x)\beta(y)\; (\!\!\mod p\;)$.
\end{itemize}
\label{G_3_id}
\end{lemma}

\begin{proof}
Since $0\cdot a=a\cdot 0=0$, it follows that $R$ is a zero-symmetric nearring if and only if $0=0\cdot b=a\alpha(0)+b\beta(0)$ or equivalently $\alpha(0)=\beta(0)=0$. Moreover, since $b=ab=a\alpha(a)+b\beta(a)$, we have $\alpha(a)=0$ and $\beta(a)=1$, so that statements (0) and (1) hold.

Further, using the left distributive law, we derive  $$xy=(xa)y_1+(xb)y_2=(ax_1+bx_2)y_1+(a\alpha(x)+b\beta(x))y_2.$$ We have also $(ax_1+bx_2)y_1=ax_1y_1+bx_2y_1$ and $(a\alpha(x)+b\beta(x))y_2=a\alpha(x)y_2+b\beta(x)y_2$. Thus $xy=a(x_1y_1+\alpha(x)y_2)+b(x_2y_1+\beta(x)y_2)$ and so statement~(2) holds.

By Corollary~\ref{cor_1} $L=\langle ap\rangle +\langle b \rangle$. Since $xL\subseteq L$ for each $x\in R$ by Lemma~\ref{prop}, we have $xb=a\alpha(x)+b\beta(x)\in L$ whence $\alpha(x)\equiv 0\; (\!\!\mod p)$ for each $x\in R$.

Finally, the associativity of multiplication in $R$ implies that
$$(xy)b=x(yb)=a\alpha(xy)+b\beta(xy).$$

Furthermore, substituting $yb=a\alpha(y)+b\beta(y)$ instead of $y$ in formula~$\rm{(***\;*)}$, we also have
$$x(yb)=a(x_1\alpha(y)+\alpha(x)\beta(y))+b(x_2\alpha(y)+\beta(x)\beta(y)).$$

Comparing the coefficients under $a$ and $b$ in two expressions obtained for $x(yb)$, we derive statements~(3) and (4) of the lemma.
\end{proof}

\begin{lemma}
$\langle ap \rangle$ is an ideal of $R$.
\label{id}
\end{lemma}

\begin{proof}
Let $x=ax_1+bx_2$, $y=ay_1+by_2\in R$, $z=apz_1\in \langle ap \rangle$. Check whether $\langle ap \rangle$ is an ideal of $R$, i.e. $(z+x)y-xy\in \langle ap \rangle$. Using formula~$\rm{(***\;*)}$,
we have
$$\begin{array}{c}
(z+x)y-xy=(apz_1+ax_1+bx_2)(ay_1+by_2)-(ax_1+bx_2)(ay_1+by_2)=\\
=(a(pz_1+x_1)+bx_2)(ay_1+by_2)-(ax_1+bx_2)(ay_1+by_2)=\\
=a((pz_1+x_1)y_1+\alpha(z+x)y_2)+b(x_2y_1+\beta(z+x)y_2)-\\
-(a(z_1y_1+\alpha(x)y_2)+b(x_2y_1+\beta(x)y_2))=\\
=a(pz_1y_1+(\alpha(z+x)-\alpha(x))y_2)+b(\beta(z+x)-\beta(x))y_2=\\
=a(pz_1y_1+(\alpha(z+x)-\alpha(x))y_2)\in \langle ap \rangle,
\end{array}$$
since $\beta(z+x)-\beta(x)=0$.

Therefore, $\langle ap \rangle$ is an ideal of $R$.
\end{proof}

\begin{lemma}
Let $R$ be a local nearring, which is not nearfield, whose additive group of $R^+$ is isomorphic to $G_3$. If $x=ax_1+bx_2$, $y=ay_1+by_2\in R$, then the mappings $\alpha\colon R\to \mathbb Z_{p^2}$ and $\beta\colon R\to \mathbb Z_{p}$ from $\rm{(***\;*)}$ can be one of the following{\rm :}
\begin{itemize}
 \item[\rm{(1)}] $\alpha(x)=0$ and $\beta(x)={x_1}^i\; (\!\!\mod p)$ ($0<i<p$){\rm ;}
 \item[\rm{(2)}] $\alpha(x)=0$ and $\beta(x)=1${\rm ;}
 \item[\rm{(3)}] $\alpha(x)=px_2$ and $\beta(x)\equiv x_1\; (\!\!\mod p)$.
\end{itemize}
\label{functions1}
\end{lemma}

\begin{proof}
1) Since zero-symmetric local nearrings of order $p^2$ are classified in~\cite{CM_68} it follows that all non-isomorphic factor-nearrings with derived subgroup ${N=R/\langle ap \rangle}$ are described. That is, you can apply the multiplication formula from the specified work, pre-adapting it for the left local nearrings. It is clear that $N^+=\langle \overline{a} \rangle + \langle \overline{b}\rangle$. Namely, let $\overline{x}=\overline{a}x_1+\overline{b}x_2$ and $\overline{y}=\overline{a}y_1+ \overline{b}y_2$ be elements of $N$, then

$$\overline{x} \overline{y}=\overline{a}x_1y_1+\overline{b}(x_2y_1+\rho(x_1)y_2).$$

Moreover, by \cite[Theorem~1.6]{CM_68} $\rho$ takes one of $p-1$ values for zero-symmetric nearrings.

On the other hand, by the formula $\rm{(***\;*)}$ we have $\overline{x} \overline{y}=\overline{a}x_1y_1+\overline{b}(x_2y_1+\beta(x)y_2)$.

Equating the coefficients for generators, we obtain: $beta(x)y_2=\rho(x_1)y_2$. Hence $\beta(x)=\rho(x_1)$. So there exist $p-1$ different zero-symmetric nearrings.

2) It is obvious that the multiplication~$\rm{(***\;*)}$ with the functions $\alpha(x)=0$ and $\beta(x)=1$ is the constant nearring multiplication.

3) By direct checking of conditions (3) and (4) of Lemma~\ref{G_3_id} for $\alpha(x)=px_2$ and $\beta(x)\equiv x_1\; (\!\!\mod p)$ the multiplication~$\rm{(***\;*)}$, we prove the lemma.
\end{proof}

So, we have examples of $p+1$ nearring multiplications and we formulate the following result.

\begin{theorem}
There exist at least $p+1$ non-isomorphic local nearrings on each metacyclic abelian groups of order $p^3$.
\label{theorem_3}
\end{theorem}

\begin{proof}
Nearrings $R=(R,+,\cdot)$ with functions $\alpha(x)=0$ and ${\beta(x)={x_1}^i\;(\!\!\mod p)}$ ($0<i<p$) are non-isomorphic zero-symmetric local nearrings according to the above and the paper~\cite{CM_68}. It is easy to check that the local nearring $R=(R,+,\cdot)$ with functions $\alpha(x)=0$ and $\beta(x)=1$ is non-zero-symmetric. It is obvious that $R=(R,+,\cdot)$ with functions $\alpha(x)=px_2$ and ${\beta(x)\equiv x_1\;(\!\!\mod p)}$ is zero-symmetric and non-isomorphic to the nearrings considered above. Therefore, there exist at least $p+1$  non-isomorphic local nearrings on each metacyclic abelian groups of order $p^3$.
\end{proof}

From the packages SONATA and LocalNR we have the following number of non-isomorphic local nearrings.
\begin{table}[h]
\centering
\begin{tabular}{|c|c|}
\hline
$StructureDescription(R^+)$ & $n(R^+)$\\
 \hline
$C_9\times C_3$                     &  13\\
 \hline
$C_{25}\times C_5$                  &  31\\
 \hline
$C_{49}\times C_7$                  &  31\\
 \hline
$C_{121}\times C_{11}$              &  47\\
   \hline
\end{tabular}
\end{table}

\section*{Acknowledgement}
This work was partially supported by the Polish Academy of Sciences (PAN) and National Academy of Sciences (NAS). The authors would like to thank the Institute of Mathematics of the Polish Academy of Sciences for hosting after evacuation from Ukraine. The authors are grateful IIE-SRF for their support of our fellowship at the University of Warsaw.

\end{document}